\overfullrule=0pt
\centerline {\bf A multiplicity theorem for a class of nonlocal elliptic equations involving even nonlinearities}\par
\bigskip
\bigskip
\centerline {BIAGIO RICCERI}\par
\bigskip
\bigskip
{\bf Abstract.} In this paper, on a bounded domain $\Omega\subset {\bf R}^n$, we consider a nonlocal problem of the type
$$\cases{-\Delta u=Q(u,\lambda)f(u)
& in $\Omega$ \cr & \cr u=0 & on $\partial\Omega$\cr}$$
where $Q:H^1_0(\Omega)\times {\bf R}\to {\bf R}$, proving, under suitable assumptions, the existence of at least three weak solutions for
each $\lambda$ running in a suitable interval.
\par
\bigskip
\bigskip
\bigskip
\bigskip
Let $\Omega\subset {\bf R}^n$ be a bounded smooth domain. On the space $H^1_0(\Omega)$ we consider the usual norm
$$\|u\|=\left(\int_{\Omega}|\nabla u(x)|^2dx\right)^{1\over 2}\ .$$
We set
$$\lambda_1:=\inf_{u\in H^1_0(\Omega)\setminus \{0\}}{{\|u\|^2}\over {\int_{\Omega}|u(x)|^2dx}}.$$
Let us introduce some notations.\par
\smallskip
If $f:\Omega\times {\bf R}\to {\bf R}$ is a Carath\'eodory function, we set
$$F(x,t)=\int_0^tf(x,s)ds$$
for all $(x,t)\in \Omega\times {\bf R}$.\par
\smallskip
If $n\geq 2$, we denote by ${\cal A}$ the family of all Carath\'eodory functions $f:\Omega\times {\bf R}\to {\bf R}$ such that
$$\sup_{(x,t)\in \Omega\times {\bf R}}{{|f(x,t)|}\over {1+|t|^q}}$$
where $q>0$, with $q<{{n+2}\over {n-2}}$ if $n>2$. If $n=1$, ${\cal A}$ denotes the family of all Carath\'eodory functions 
$f:\Omega\times {\bf R}\to {\bf R}$ such that $\sup_{|t|\leq r}|f(\cdot,t)|\in L^1(\Omega)$ for all $r>0$.\par
\smallskip
Given a continuous function $Q:{\bf R}\to {\bf R}$, and two functions $f, g\in {\cal A}$, we are interested in the Dirichlet problem
$$\cases{-\Delta u=Q\left(\int_{\Omega}F(x,u(x))dx\right)f(x,u)+g(x,u) & in $\Omega$ \cr & \cr u=0 & on $\partial\Omega$.\cr}\eqno{(P)}$$
\smallskip
As usual, a weak solution of $(P)$ is any $u\in H^1_0(\Omega)$
such that
$$\int_{\Omega}\nabla u(x)\nabla v(x)dx-Q\left(\int_{\Omega}F(x,u(x))dx\right)\int_{\Omega}f(t,u(x))v(x)dx-\int_{\Omega}g(t,u(x))v(x)dx=0$$
for all $v\in H^1_0(\Omega)$.
\smallskip
This is a nonlocal problem due to the presence of $Q$ (we refer to [1] for an introduction to the subject).
\smallskip
Let $J, \Phi:H^1_0(\Omega)\to {\bf R}$ be the functionals defined by
$$J(u)=\tilde Q\left(\int_{\Omega}F(x,u(x))dx\right)+\int_{\Omega}G(x,u(x))dx\ ,$$
$$\Phi(u)={{1}\over {2}}\|u\|^2-J(u),$$
where $\tilde Q(t)=\int_0^tQ(s)ds$.\par
\smallskip
Since $f, g\in {\cal A}$, by classical results, the functional $J$ is $C^1$ and sequentially weakly continuous, and the weak solutions of $(P)$
agree with the critical points of $\Phi$ in $H^1_0(\Omega)$.\par
\medskip
The aim of this very short note is to prove the following result:\par
\medskip
THEOREM 1. {\it Let $\alpha:{\bf R}\to {\bf R}$ be a derivable strictly convex function such that $\lim_{|\lambda|\to +\infty}
\alpha(\lambda)=+\infty$, $\alpha(0)>0$, $\alpha'(0)=0$, $\alpha'({\bf R})={\bf R}$.
Moreover, let $f\in {\cal A}$ and let $\gamma:{\bf R}\to {\bf R}$ be a continuous, odd
function, vanishing only at $0$, with $\gamma(1)>0$, such that\par
\noindent
$(a)$\hskip 5pt for almost every $x\in \Omega$, the function $f(x,\cdot)$ is even in ${\bf R}$;\par
\noindent
$(b)$\hskip 5pt $\hbox {\rm meas}(\{x\in \Omega : \sup_{t\in {\bf R}}|f(x,t)|>0\})>0$;\par
\noindent
$(c)$\hskip 5pt $\lim_{|t|\to +\infty}{{\sup_{x\in \Omega}|F(x,t)|}\over {t^2}}=0$;\par
\noindent
$(d)$\hskip 5pt $\limsup_{\to +\infty}{{\int_0^t\gamma(s)ds}\over {t}}<+\infty$.\par
Then, for each $\mu>0$ large enough, exists $\epsilon_{\mu}^*>0$ with the following property:
for each $g\in {\cal A}$ satisfying
$$\int_{\Omega}\sup_{t\in {\bf R}}g(x,t)dx-\int_{\Omega}\inf_{t\in {\bf R}}g(x,t)dx<\epsilon_{\mu}^*\ ,$$
 there exists an open interval $A\subseteq {\bf R}$ and a number $\rho>0$ such that, for each $\lambda\in A$, the problem
$$\cases{-\Delta u=\left(\lambda+\mu\alpha(\lambda)\gamma\left(\int_{\Omega}F(x,u(x))dx\right)\right)f(x,u)+g(x,u) & in $\Omega$ \cr & \cr u=0 & on $\partial\Omega$\cr}$$
has at least three weak solutions whose norms in $H^1_0(\Omega)$ are less than $\rho$.}\par
\smallskip
PROOF. First, we show that $\int_{\Omega}F(x,u(x))>0$ for some $u\in H^1_0(\Omega)$. To this end,
set
$$D=\left \{x\in \Omega : \sup_{\xi\in {\bf R}}F(x,\xi)>0\right\}\ .$$
By $(b)$, meas$(D)>0$. Then, by  the Scorza-Dragoni theorem,
there exists a compact set $K\subset D$, of positive measure, such that the restriction
of $F$ to $K\times {\bf R}$ is continuous. Fix a point $\hat x\in K$ such that
the intersection of $K$ and any ball centered at $\hat x$ has a positive measure.
Choose $\hat \xi\in {\bf R}$ so that $F(\hat x,\hat\xi)>0$. By continuity, there
is $r>0$ such that
$$F(x,\hat\xi)>0$$
for all $x\in K\cap B(\hat x,r)$. Set
$$M=\sup_{(x,\xi)\in \Omega\times [-|\hat \xi|,|\hat \xi|]}|F(x,\xi)|\ .$$
Since $f\in {\cal A}$, we have $M<+\infty$. Next, choose an open set $\tilde\Omega$
such that
$$K\cap B(\hat x,r)\subset\tilde\Omega\subset\Omega$$
and
$$\hbox {\rm meas}(\tilde\Omega\setminus (K\cap B(\hat x,r))<{{\int_{K\cap B(\hat x,r)}F(x,\hat \xi)dx}\over {M}}\ .$$
Finally, choose a function $\tilde u\in H^1_0(\Omega)$ such that
$$\tilde u(x)=\hat\xi$$
for all $x\in K\cap B(x,r)$,
$$\tilde u(x)=0$$
for all $x\in \Omega\setminus\tilde\Omega$ and
$$|\tilde u(x)|\leq |\hat\xi|$$
for all $x\in\Omega$. Thus, we have
$$\int_{\Omega}F(x,\tilde u(x))dx=\int_{K\cap B(\hat x,r)}F(x,\hat \xi)dx+\int_{\tilde\Omega\setminus (K\cap B(\hat x,r)}F(x,\tilde u(x))dx$$
$$>\int_{K\cap B(\hat x,r)}F(x,\hat \xi)dx-M\hbox {\rm meas}(\tilde\Omega\setminus (K\cap B(\hat x,r))>0\ .$$
So, since $\Gamma(t):=\int_0^t\gamma(s)ds>0$ if $t\neq 0$, the number
$$\eta^*:=\sup_{u\in H^1_0(\Omega)}{{\Gamma\left(\int_{\Omega}F(x,u(x))dx\right)}\over {\|u\|^2}}$$
is positive. 
Fix $\mu$ so that
$$\mu>{{1}\over {2\alpha(0)\eta^*}}.\eqno{(1)}$$
Next, fix $a>0$ and consider the functional $P:H^1_0(\Omega)\times {\bf R}\to {\bf R}$ defined by
$$P(u,\lambda):={{1}\over {2}}\|u\|^2-\lambda\int_{\Omega}F(x,u(x))dx-\mu\alpha(\lambda)
\left(\Gamma\left(\int_{\Omega}F(x,u(x))dx\right)+a\right).$$
We claim that
$$\sup_{{\bf R}}\inf_{H^1_0(\Omega)}P<\inf_{H^1_0(\Omega)}\sup_{{\bf R}}P.\eqno{(2)}$$
Arguing by contradiction, suppose that
$$\sup_{{\bf R}}\inf_{H^1_0(\Omega)}P=\inf_{H^1_0(\Omega)}\sup_{{\bf R}}P.\eqno{(3)}$$
Since 
$$P(0,\lambda)=-a\mu\alpha(\lambda), $$
we have
$$\lim_{|\lambda|\to +\infty}\inf_{u\in H^1_0(\Omega)}P(u,\lambda)=-\infty.$$
But the function $\lambda\to \inf_{u\in H^1_0(\Omega)}P(u,\lambda)$ is upper semicontinuous and hence
there exists
$\lambda^*\in {\bf R}$ such that
$$\sup_{\lambda\in {\bf R}}\inf_{u\in H^1_0(\Omega)}P(u,\lambda)=\inf_{u\in H^1_0(\Omega)}P(u,\lambda^*).\eqno{(4)}$$
On the othe hand, the functional $u\to \sup_{\lambda\in {\bf R}}P(u,\lambda)$ is sequentially weakly lower semicontinuous. Let us
check that it is also coercive. Fix $\lambda\in {\bf R}$. By $(d)$, there are two numbers $c_1, c_2>0$ such that
$$\Gamma(t)\leq c_1t+c_2$$
for all $t\geq 0$. Now, fix $\epsilon >0$ so that
$$\epsilon<{{\lambda_1}\over {2(|\lambda|+\mu\alpha(\lambda)c_1)}}.$$
By $(c)$, taking into account that $f\in {\cal A}$, there is another number $c_3>0$ such that
$$|F(x,t)|\leq \epsilon t^2+c_3$$
for all $t\in {\bf R}$. Then, since $\Gamma$ is even and increasing in $[0,+\infty[$, for suitable constants $c_4, c_5>0$,
for each $u\in H^1_0(\Omega)$, we have
$$|\lambda|\left|\int_{\Omega}F(x,u(x))dx\right|+\mu|\beta(\lambda)\Gamma\left(\int_{\Omega}F(x,u(x))dx\right)\leq
|\lambda|\int_{\Omega}|F(x,u(x))|dx+\mu\alpha(\lambda)\Gamma\left(\int_{\Omega}|F(x,u(x))|dx\right)$$
$$\leq (|\lambda|+\mu\alpha(\lambda)c_1)\int_{\Omega}|F(x,u(x))|dx+c_4\leq (|\lambda|+\mu\alpha(\lambda)c_1)\epsilon\int_{\Omega}|u(x)|^2dx+c_5$$
$$\leq {{(|\lambda|+\mu\alpha(\lambda)c_1)\epsilon}\over {\lambda_1}}\|u\|^2+c_5.$$
Consequently, for another constant $c_6>0$, we have
$$P(u,\lambda)\geq \left({{1}\over {2}}-{{(|\lambda|+\mu\alpha(\lambda)c_1)\epsilon}\over {\lambda_1}}\right)\|u\|^2-c_6.$$
Hence, since the factor multiplying $\|u\|^2$ is positive, we have
$$\lim_{\|u\|\to +\infty}P(u,\lambda)=+\infty,$$
as claimed. Jointly with the sequential weak lower semicontinuity, this implies that
there exists $u^*\in H^1_0(\Omega)$ such that
$$\inf_{u\in H^1_0(\Omega)}\sup_{\lambda\in {\bf R}}P(u,\lambda)=\sup_{\lambda\in {\bf R}}P(u^*,\lambda).\eqno{(5)}$$
So, from $(4), (5), (6)$, we obtain
$$P(u^*,\lambda^*)=\inf_{u\in H^1_0(\Omega)}P(u,\lambda^*)=\sup_{\lambda\in {\bf R}}P(u^*,\lambda).\eqno{(6)}$$
Now, set $g:=-\alpha'$. By assumption, $g$ is strictly decreasing and surjective. From this and $(5)$ we readily
infer that
$$\lambda^*=g^{-1}\left({{\int_{\Omega}F(x,u^*(x))dx}\over {\mu(\Gamma\left(\int_{\Omega}F(x,u^*(x))dx\right)+a)}}\right).$$
Set
$$k^*:=-\alpha\left(g^{-1}\left({{\int_{\Omega}F(x,u^*(x))dx}\over {\mu(\Gamma\left(\int_{\Omega}F(x,u^*(x))dx\right)+a)}}\right)\right).$$
From $(6)$ again, we have
$${{1}\over {2}}\|u^*\|^2-\lambda^*\int_{\Omega}F(x,u^*(x))dx+\mu k^*\left(\Gamma\left(\int_{\Omega}F(x,u^*(x))dx\right)+a\right)$$
$$=\inf_{u\in H^1_0(\Omega)}\left({{1}\over {2}}\|u\|^2-\lambda^*\int_{\Omega}F(x,u(x))dx+\mu k^*\left(\Gamma\left(\int_{\Omega}F(x,u(x))dx\right)+a\right)\right).\eqno{(7)}$$
Now, suppose that $\int_{\Omega}F(x,u^*(x))dx\neq 0$. Notice that, if $\int_{\Omega}F(x,u^*(x))dx>0$, since $g^{-1}$ is strictly decreasing and vanishes at $0$,  we have
$\lambda^*<0$. For the same reasons, if $\int_{\Omega}F(x,u^*(x))dx<0$ then $\lambda^*>0$. That is
to say, $\lambda^*\int_{\Omega}F(x,u^*(x))dx<0$. Then, since $F(x,\cdot)$ is odd  and $\Gamma$ is even, if we choose $v=-u^*$, we have
$${{1}\over {2}}\|v\|^2-\lambda^*\int_{\Omega}F(x,v(x))dx+\mu k^*\left(\Gamma\left(\int_{\Omega}F(x,v(x))dx\right)+a\right)$$
$$={{1}\over {2}}\|u^*\|^2+\lambda^*\int_{\Omega}F(x,u^*(x))dx+\mu k^*\left(\Gamma\left(\int_{\Omega}F(x,u^*(x))dx\right)+a\right)$$
$$<{{1}\over {2}}\|u^*\|^2-\lambda^*\int_{\Omega}F(x,u^*(x))dx+\mu k^*\left(\Gamma\left(\int_{\Omega}F(x,u^*(x))dx\right)+a\right).\eqno{(8)}$$
From $(7)$ and $(8)$, it follows that  $\int_{\Omega}F(x,u^*(x))dx=0$. Assume that $u^*\neq 0$.  We then have
$$-\mu\alpha(0)a<{{1}\over {2}}\|u^*\|-\mu\alpha(0)a.\eqno{(9)}$$
Then, observing that $\lambda^*=0$ and $k^*=-\alpha(0)$, from $(7)$ and $(9)$ we infer that $u^*=0$. In view of $(1)$,
there is some $w\in H^1_0(\Omega)$ so that
$$-{{1}\over {2\mu\beta(0)}}<{{\Gamma\left(\int_{\Omega}F(x,w(x))dx\right)}\over {\|w\|^2}}.$$
From this, we get
$${{1}\over {2}}\|w\|^2-\mu\alpha(0)\left(\Gamma\left(\int_{\Omega}F(x,w(x)\right)dx+a\right)<-\mu\alpha(0)a. \eqno{(10)}$$
But $(10)$ contradicts $(7)$. This contradiction shows that $(2)$ is true. Now, put
$$\epsilon_{\mu}^*:=\inf_{H^1_0(\Omega)}\sup_{{\bf R}}P-\sup_{{\bf R}}\inf_{H^1_0(\Omega)}P.$$
Fix any $g\in {\cal A}$ such that
$$\int_{\Omega}\sup_{t\in {\bf R}}G(x,t)dx-\int_{\Omega}\inf_{t\in {\bf R}}G(x,t)dx<\epsilon_{\mu}^*.$$
Consider the functional $\varphi:H^1_0(\Omega)\to {\bf R}$ defined by
$$\varphi(u):=-\int_{\Omega}G(x,u(x))dx.$$
Since
$$\int_{\Omega}\inf_{t\in {\bf R}}G(x,t)dx\leq \varphi(u)\leq \int_{\Omega}\sup_{t\in {\bf R}}G(x,t)dx,$$
we have
$$\sup_{H^1(\Omega)}\varphi-\inf_{H^1(\Omega)}\varphi\leq \int_{\Omega}\sup_{t\in {\bf R}}G(x,t)dx-\int_{\Omega}\inf_{t\in {\bf R}}G(x,t)dx<
\epsilon_{\mu}^*.$$
By [3], [4], we have the strict inequality
$$\sup_{{\bf R}}\inf_{H^1_0(\Omega)}(P-\varphi)<\inf_{H^1_0(\Omega)}\sup_{{\bf R}}(P-\varphi).$$
Then, since $P(u,\cdot)$ is continuous and concave, and $P(\cdot,\lambda)-\varphi(\cdot)$ is weakly lower semicontinuous, coercive and satisfies the Palais-Smale condition (see Example 38.25 of [5]), we can apply Theorem 3 of [2] from which the conclusion follows directly.\hfill $\bigtriangleup$
\medskip
REMARK 1. - We are not aware of known results close enough to Theorem 1 so that a proper comparison can be made.
\bigskip
\bigskip
{\bf Acknowledgements.} 
The author has been supported by the Gruppo Nazionale per l'Analisi Matematica, la Probabilit\`a e 
le loro Applicazioni (GNAMPA) of the Istituto Nazionale di Alta Matematica (INdAM).

\bigskip
\bigskip
\bigskip
\centerline {\bf References}\par
\bigskip
\bigskip
\noindent
[1]\hskip 5pt M. CHIPOT, {\it Remarks on some class of nonlocal elliptic problems}, Recent Advances on Elliptic and Parabolic Issues, World Scientific, Singapore, 2006, 79-102.\par
\smallskip
\noindent
[2]\hskip 5pt  B. RICCERI, {\it On a three critical points theorem}, Arch. Math. (Basel), {\bf 75} (2000), 220-226.\par
\smallskip
\noindent
[3]\hskip 5pt B. RICCERI, {\it A more complete version of a minimax theorem}, Appl. Anal. Optim., {\bf 5} (2021), 251-261.\par
\smallskip
\noindent
[4]\hskip 5pt B. RICCERI, {\it Addendum to  ``A more complete version of a minimax theorem"}, Appl. Anal. Optim., {\bf 6} (2022), 195-197.\par
\smallskip
\noindent
[5]\hskip 5pt E. ZEIDLER, {\it Nonlinear functional analysis and its applications}, vol. III, Springer-Verlag, 1985.\par
\bigskip
\bigskip
\bigskip
\bigskip
Department of Mathematics and Informatics\par
University of Catania\par
Viale A. Doria 6\par
95125 Catania, Italy\par
{\it e-mail address}: ricceri@dmi.unict.it

\bye